\def\E{\end{document}}
\documentclass[12pt]{article}
\usepackage{amssymb}
\usepackage{amssymb,amsmath,color}
\usepackage{cases}
\usepackage{mathrsfs}
\usepackage{floatrow}
\usepackage{graphicx}

\newtheorem{theorem}{Theorem}[section]
\newtheorem{definition}{Definition}[section]
\newtheorem{lemma}{Lemma}[section]

\newtheorem{corollary}{Corollary}[section]
\newtheorem{remark}{Remark}[section]

\catcode`@=11 \@addtoreset{equation}{section} \catcode`@=12

\begin{document}

\title{\bf The stability and rapid exponential stabilization of heat equation in non-cylindrical  domain
\thanks{This work is partially
supported by the NSF of China under grants 11471070 and 11371084. }}
\author{ Lingfei Li\thanks{School of Mathematics and Statistics, Northeast Normal
University, Changchun 130024, China; School of Science, Northeast Electric Power
University, Jilin 132012, China. E-mail address: lilf320@163.com.
},\quad Yujing Tang\thanks{Experimental High School of Qiqihar City, Qiqihar 161000, China. E-mail address:
820299877@qq.com  }\quad  and\quad Hang Gao
\thanks{Corresponding author, School of Mathematics and Statistics, Northeast Normal
University, Changchun 130024, China. E-mail address: hanggao2013@126.com}}
\date{}
\maketitle

\begin{abstract}
 This paper is devoted to the study of the stability and stabilizability of heat equation in non-cylindrical domain. The interesting thing is that there is a class of initial values such that the system is no longer exponentially stable. The system is only polynomially stable or only analogously exponentially stable. Then, the rapid exponential stabilization of the system is obtained by the backstepping method.
\end{abstract}


\section{Introduction and Main results}
It is well known that the stability and stabilization problems for both linear and nonlinear partial differential equations have been studied extensively(see [1-4],  references therein). In general, the study of stability for a given system is along the following way. We first concern whether the solution is stable or not. Then, if it is unstable, we try to find a control to stabilize the system. And if the solution decays in a slower rate, one wants to force the solution to decay with arbitrarily prescribed decay rates. They are called stability, stabilization and rapid stabilization, respectively.

There are many methods to study the stabilization including pole placement, the control Lyapunov function method, the backstepping method and so on.  The backstepping method which was initiated in [5] and [6] has been used as a standard method for finite dimensional control systems. The application of continuous backstepping method to parabolic equations was first given in [7] and [8]. In the last decade, the backstepping method has been widely used to study the stability of partial differential equations, such as  wave equation, Korteweg-de Vries equation, and Kuramoto-Sivashinsky equation (see [9]-[11]). The scheme of backstepping method is as follows. Initially, an exponentially stable target system is established. Then, the PDE system is transformed to the target system by using a Volterra transformation. At the same time, the PDE describing the transformation kernel is obtained. Thus the stabilization problem is converted to the well-posedness problem of the kernel and the invertibility of transformation. Finally, An explicit full state feedback control is given to demonstrate successful stabilization of the unstable system.

Many problems in the real world involve non-cylindrical regions such as controlled annealing of a solid in a fluid medium([12]), vibration control of an extendible flexible beam([12]), phase change and heat transfer. The PDE in the non-cylindrical domain has close relation with equations of time-dependent coefficients. The typical system with time-dependent coefficients is given by the Czochralski crystal growth problem ([13]), which is presented as a heat equation with time-dependent coefficients. Hence, the research on the problem of non-cylindrical region has important practical significance. Time-dependent domain and coefficients lead to more complexities and difficulties. In most papers concerning stabilization, systems are considered in cylindrical domains. Only few results have been found for systems in non-cylindrical domains(see [12][14][15]). In [12], the authors considered the problem of the stabilization and the control of distributed systems with time-dependent spatial domains. The evolution of the spatial domains with time is described by a finite-dimensional system of ordinary differential equations depending on the control. Namely, the dynamical behavior of the distributed system is controlled by manipulating its spatial domain. The length of the time interval is finite. In [14], stabilization of heat equation with time-varying spatial domain is discussed. The well-posedness of the kernel which depends on time $t$ is proved by successive approximation. There is a restriction on the boundary of moving domain. More precisely, the boundary function $l(t)$ is analytic and the $j$th partial derivative satisfies:
 $$ |\partial_{t}^{j}l(t)|\leq M^{j+1}j!, j\geq0.$$
If the boundary is unbounded with respect to $t$, then the assumption in [14] does not hold.
In [15], the authors extended the backstepping-based observer design in [14] to the state estimation of parabolic PDEs  with time-dependent spatial domain.
 As far as we know, the stability and stabilizability of parabolic equation in unbounded non-cylindrical domain have not been discussed in detail yet.

In this paper, we mainly focus on the stability and rapid stabilization of the one-dimensional heat equation in non-cylindrical domain. Let us begin with the following system
\begin{equation}
\label{b1}
\left\{\begin{array}{ll}
 u_{t}-u_{xx}=0
&\mbox{in}\ {Q}_t,\\[3mm]
u(0,t)=0,u(l(t),t)=0, &\mbox{in}\ (0,+\infty), \\[3mm]
u(x,0)=u_{0}(x), &\mbox{in}\ \Omega=(0,1),
\end{array}
\right.
\end{equation}
where ${Q}_t=\{(x,t)| x\in (0,l(t))=\Omega_{t},t\in(0,+\infty), l(t)=(1+kt)^{\alpha},\alpha>0, k>0\}$.
The existence and uniqueness of solutions to the parabolic equations in non-cylindrical domains are investigated in [16]. The following Lemma can be proved by the method in [16].
\begin{lemma}
\label{b3}
If $u_{0}(x)\in L^2(\Omega)$, then system (1.1) has a unique weak solution $u$ in the following space $$ C([0,+\infty);L^2(\Omega_{t}))\cap  L^2(0,+\infty;H^1_{0}(\Omega_{t})),$$and there exists a positive constant $C$ independent of $u_{0}$ and ${Q}_t$ such that
$$\|u\|_{L^{\infty}(0,+\infty;L^2(\Omega_{t}))}+ \|\nabla u\|_{L^{2}({Q}_t)}\leq C\|u_{0}\|_{L^2(\Omega)}.$$
Moreover, if  $u_{0}(\cdot)\in H^1_{0}(\Omega)$, then system (1.1) has a unique strong solution $u$ in the class
$$ C([0,+\infty);H^1_{0}(\Omega_{t}))\cap  L^2(0,+\infty;H^2\cap H^1_{0}(\Omega_{t}))\cap H^1(0,+\infty;L^2(\Omega_{t})),$$
and there exists a constant independent of $u_{0}$ and ${Q}_t$ such that
$$\|u\|_{L^{\infty}(0,+\infty;H^1_{0}(\Omega_{t}))}+ \|u_t\|_{L^{2}({Q}_t)}\leq C\|u_{0}\|_{H^1_{0}(\Omega)}.$$
\end{lemma}

The proof of Lemma 1.1 will be given in Appendix. In what follows, the definitions of the stability for system (1.1) are given.
\begin{definition}\label{0}
System $(1.1)$ is said to be exponentially stable, if
for any given $u_{0}\in L^2(\Omega)$, there exist  $C>0$,  $ \alpha>0$ and  $t_{0}>0$ such that for every $t\geq t_{0}$, $\|u\|_{L^{2}(\Omega_{t})}\leq C e^{-\alpha t}$. System $(1.1)$ is said to be analogously exponentially stable, if
for any given $u_{0}\in L^2(\Omega)$, there exist $C>0$, $ C_{1}>0$, $t_{0}>0$  and $\beta$ with $ 0<\beta<1$ such that for every $t\geq t_{0}$, $\|u\|_{L^{2}(\Omega_{t})}\leq C e^{-C_{1}{t}^{\beta}}$. And system $(1.1)$ is said to be polynomially stable, if for any given $u_{0}\in L^2(\Omega)$, there exist $\gamma>0$, $t_{0}>0$ and $C>0$ such that for every $t\geq t_{0}$, $\|u\|_{L^{2}(\Omega_{t})}\leq C{(\varphi(t))^{-\gamma}},$ where $\varphi(t)$ is a polynomial with respect to $t$ .
\end{definition}

It is easy to see that the exponentially stable system must be analogously exponentially stable and the analogously exponentially stable system must be polynomially stable. We can get by energy estimate that system (1.1) is polynomially (or analogously exponentially) stable. What we want to know is that whether the polynomially (or analogously exponentially) stable system is exponentially stable or not. The conclusion is that the corresponding solution ${u}$ is only polynomially stable for $\alpha\geq\frac{1}{2}$, more precisely, the $L^{2}$ norm of solution has a polynomial lower bound for some initial values. And system $(1.1)$ is only analogously exponentially stable for $0<\alpha<\frac{1}{2}$ because we can also find a lower bound of the solution. Hence, we have the following Theorems.

\begin{theorem}\label{02} System $(1.1)$ is only polynomially stable for $\alpha\geq\frac{1}{2}$.
\end{theorem}
\begin{theorem}\label{01}
System $(1.1)$ is only analogously exponentially stable for $0<\alpha<\frac{1}{2}$.
\end{theorem}

\begin{remark}\label{01}
It is well known that the heat equation on the cylindrical domain is exponentially stable. If $\alpha=1$, then the boundary is a line. It can be seen by Theorem 1.1 that the system is polynomially stable as long as the boundary of the domain is tilted a little bit, namely $k$ is small enough. When the boundary of the domain is inclined steeply so that it is close to the x-axis, i.e.,k is large enough, the system is still polynomially stable. But if the boundary is $(1+kt)^\alpha(0<\alpha<\frac{1}{2})$, system (1.1) is only analogously exponentially stable. The relation between the stability and the boundary curve is shown in Figure 1.
\end{remark}

\begin{figure}[!ht]
\centering
\includegraphics[width=0.47\textwidth,height=5.5cm]{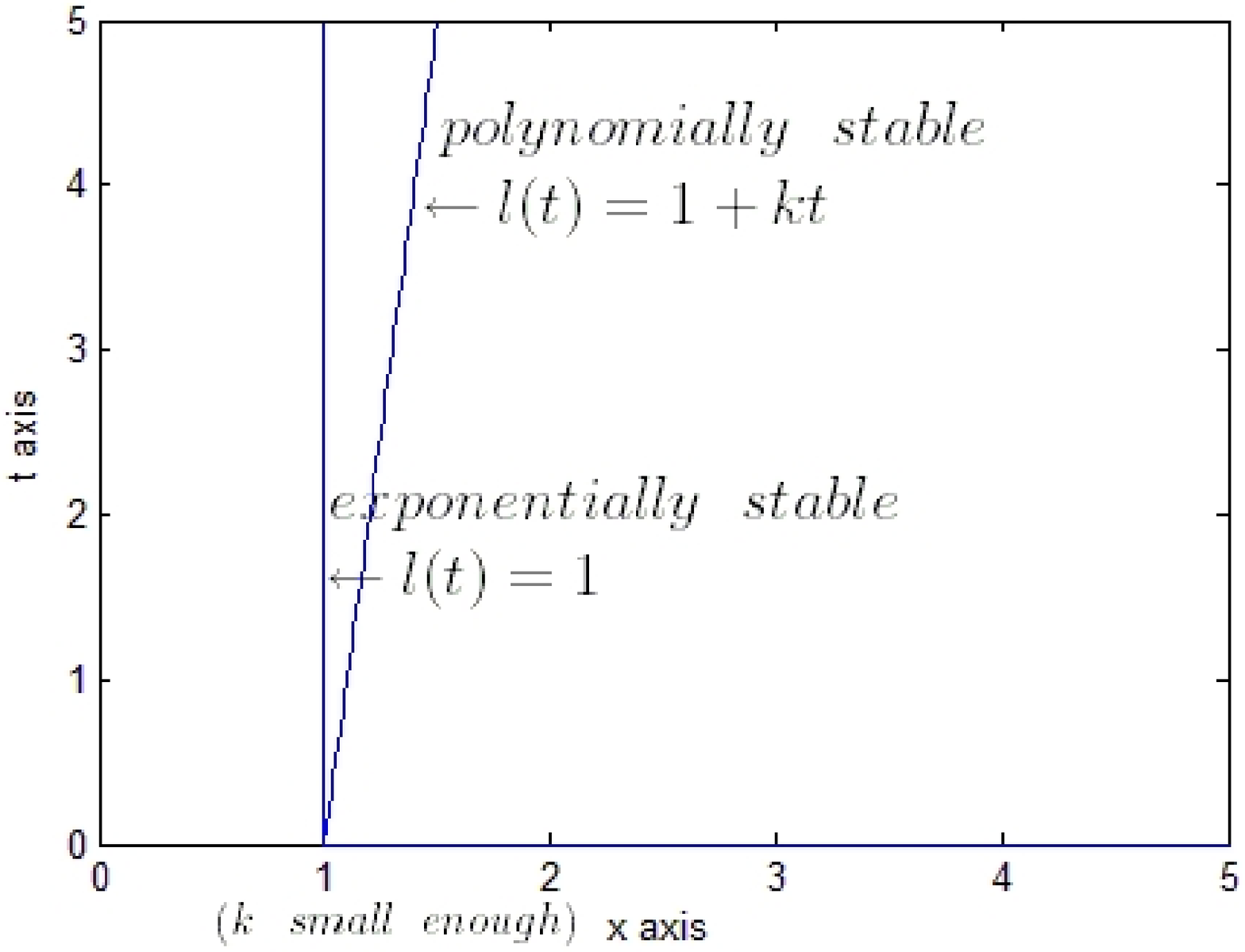}
\includegraphics[width=0.47\textwidth,height=5.5cm]{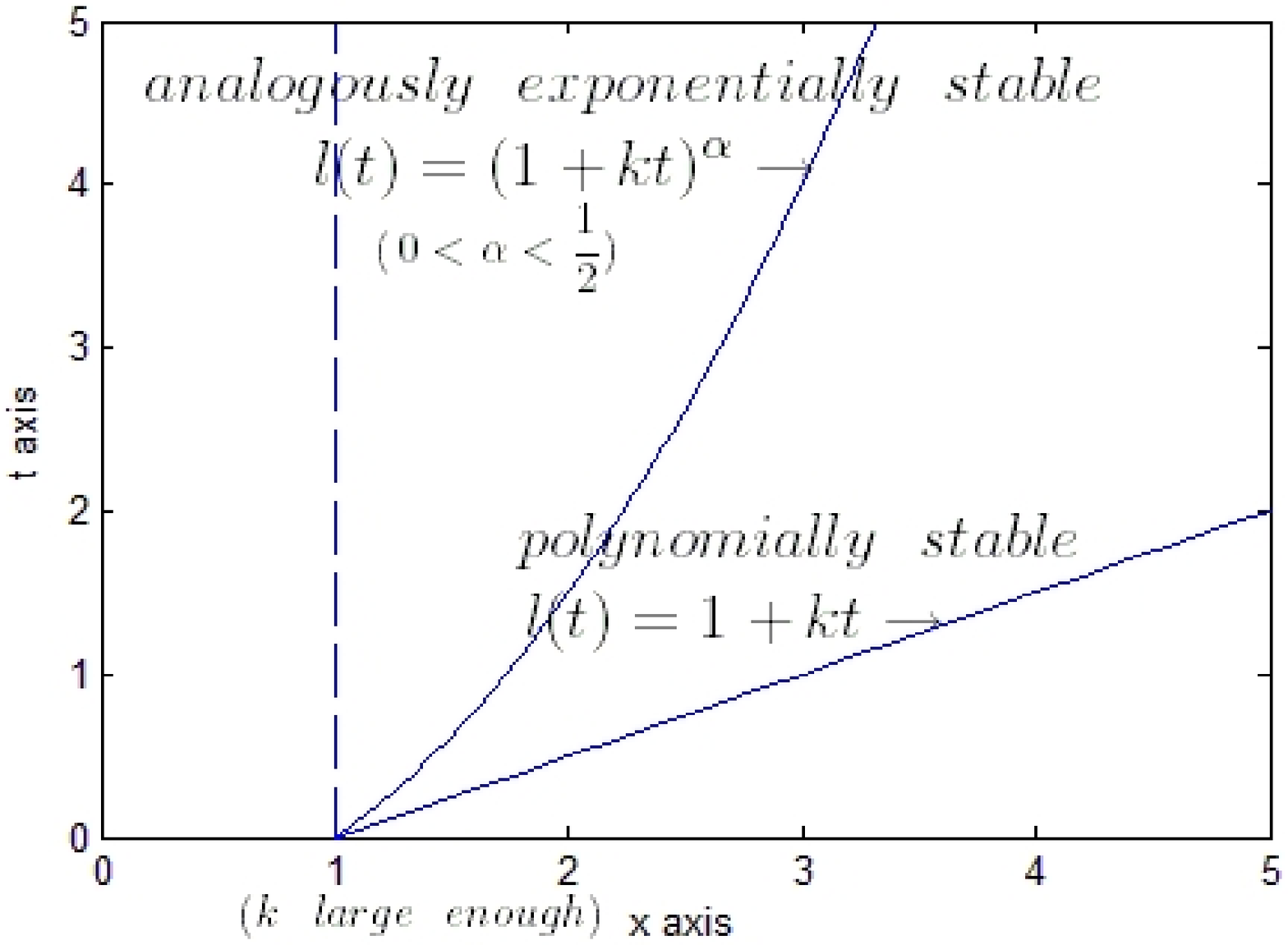}
\caption{The relation between the stability and the boundary curve}
\end{figure}

Since system (1.1) decays slower than exponential decay rate, the goal of this paper is to construct a control to force the solution to decay at a desired rate. It is easy to see that there exists an internal feedback control leading to exponential decay. However, we are interested in looking for boundary feedback control to stabilize the system exponentially. Let us consider the system with boundary control

\begin{equation}
\label{b1}
\left\{\begin{array}{ll}
 u_{t}-u_{xx}=0,
&\mbox{in}\ {Q}_t,\\[3mm]
u(0,t)=0,u(l(t),t)=U(t), &\mbox{in}\ (0,+\infty), \\[3mm]
u(x,0)=u_{0}(x), &\mbox{in}\ \Omega=(0,1),
\end{array}
\right.
\end{equation}
where ${Q}_t=\{(x,t)| x\in (0,l(t))=\Omega_{t},t\in(0,+\infty), l(t)=(1+kt)^{\alpha},0<\alpha\leq1, k>0\}$.
The well-posedness of solutions to (1.2) can be obtained from the results in [16] and Lemma 1.1.
\begin{lemma}
\label{b3}
Assume $U(t)\in H^1(0,+\infty)$ and $\sqrt{l(t)}U'(t)\in L^2(0,+\infty)$. If $u_{0}(x)\in L^2(\Omega)$, system (1.2) has a unique weak solution $u$ in the class $$ C([0,+\infty);L^2(\Omega_{t}))\cap  L^2(0,+\infty;H^1(\Omega_{t})).$$
 Moreover, if  $u_{0}(x)\in H^1(\Omega)$ and $u_{0}(1)=U(0)$, then system (1.2) has a unique strong solution $u$ in the class $$ C([0,+\infty);H^1(\Omega_{t}))\cap  L^2(0,+\infty;H^2(\Omega_{t}))\cap H^1(0,+\infty;L^2(\Omega_{t})).$$
\end{lemma}

The proof of Lemma 1.2 will be shown in Appendix.
\begin{theorem}\label{02}
Assume that $u_{0}\in L^2(\Omega),0<\alpha\leq1$. Then, there exists a boundary feedback control $U(t)\in H^{1}(0,+\infty)$ such that system $(1.1)$ decays exponentially to zero.
\end{theorem}
\begin{remark}\label{01}
When the growth order of $l(t)$ is more than 1, the linear feedback control does not work. However, we can not deduce that the system is not stabilized exponentially  by feedback control because of the variety of feedback control.
\end{remark}

This paper is organized as follows. In section 2, the stability of (1.1) is established. While in Section 3, the feedback stabilization of (1.2) for  $0<\alpha\leq 1$ is proved. At last, the appendix is given in Section 4.
\section{The stability of (1.1)}
\begin{lemma}\label{01}
Let us assume that $u_{0}\in C(\overline{\Omega})$, $u_{0}>0$,$\alpha\geq\frac{1}{2}$. Then
there is some $t_{0}$ such that the solution of system (1.1) satisfies
$$\|u\|_{L^2(0,l(t))}\geq C (1+kt)^{-\beta} ,t\geq t_{0},$$
where $ C$ depending on $k,\alpha, t_{0}$, and $ \beta>0$ depending on $k,\alpha$.
\end{lemma}

\noindent{\bf Proof of Lemma 2.1}
(1)Let $u_{0}(x)=\sin(\pi x)e^{-\frac{\alpha x^{2}}{4}}$. Assume  that $u(x,t)=\sin\frac{\pi x}{(1+kt)^{\alpha}}e^{c(x,t)}$ is the solution of the following system,
\begin{equation}
\label{b1}
\left\{\begin{array}{ll}
 u_{t}-u_{xx}=0
&\mbox{in}\ {Q}_t,\\[3mm]
u(0,t)=0,u(l(t),t)=0, &\mbox{in}\ (0,+\infty), \\[3mm]
u(x,0)=u_{0}(x), &\mbox{in}\ \Omega=(0,1),
\end{array}
\right.
\end{equation}
where $\alpha\geq\frac{1}{2}$.

By straightforward computations, we get that the partial derivatives are
\begin{equation}
\begin{array}{ll}
\frac{\partial u}{\partial t}=-\frac{\alpha \pi kx}{(1+kt)^{\alpha+1}}\cos\frac{\pi x}{(1+kt)^{\alpha}}e^{c(x,t)}+\sin\frac{\pi x}{(1+kt)^{\alpha}}e^{c(x,t)} \frac{\partial c }{\partial t},\\[3mm]

\frac{\partial u}{\partial x}=\frac{\pi}{(1+kt)^{\alpha}}\cos\frac{\pi x}{(1+kt)^{\alpha}}e^{c(x,t)}+\sin{\frac{\pi x}{(1+kt)^\alpha}}e^{c(x,t)} \frac{\partial c}{\partial x},\\[3mm]

\frac{\partial^{2} u}{\partial x^{2}}=
-\frac{\pi^{2}}{(1+kt)^{2\alpha}}\sin\frac{\pi x}{(1+kt)^{\alpha}}e^{c(x,t)}+\frac{2 \pi }{(1+kt)^{\alpha}}\cos\frac{\pi x}{(1+kt)^{\alpha}}e^{c(x,t)} \frac{\partial c}{\partial x}\\[3mm]
+\sin\frac{\pi x}{(1+kt)^{\alpha}}e^{c(x,t)}( (\frac{\partial c}{\partial x})^{2}+\frac{\partial^{2} c}{\partial x^{2}}).
\end{array}
\end{equation}
Substituting (2.2) into the equation in (2.1) and comparing the coefficients of sine and cosine, one has the equation satisfied by $c(x,t)$ as follows
\begin{equation}
\label{b6}
\left\{\begin{array}{ll}
 \frac{\partial c}{\partial x}=-\frac{\alpha k x}{2(1+kt)},\\[3mm]
\frac{\partial c}{\partial t}=-\frac{\pi^{2}}{(1+kt)^{2\alpha}}+\frac{\alpha^{2}k^{2}x^{2}}{4(1+kt)^2}-\frac{\alpha k }{2(1+kt)}.
\end{array}
\right.
\end{equation}

When $\alpha=\frac{1}{2}$, integrating (2.3) respect for $x$ and $t$ respectively, we obtain
\begin{equation}
\label{b7}
\left\{\begin{array}{ll}
  c(x,t)- c(0,t)=-\frac{\alpha kx^{2}}{4(1+kt)},\\[3mm]
 c(x,t)- c(x,0)=-\frac{\pi^{2}\ln(1+kt)}{k}-\frac{k\alpha^{2}x^{2}}{4(1+kt)}+
 \frac{k\alpha^{2}x^{2}}{4}-\frac{\alpha}{2}\ln(1+kt).
\end{array}
\right.
\end{equation}
Taking $c(0,0)=0$, then we have
$$c(x,t)=-\frac{\pi^{2}\ln(1+kt)}{k}-\frac{\alpha}{2}\ln(1+kt)-
\frac{k\alpha^{2}x^{2}}{4(1+kt)}+\frac{k\alpha^{2}x^{2}}{4}-\frac{k\alpha x^{2}}{4}.$$
Hence, the expression of the solution to (2.1) is
\begin{equation}
\label{b8}
\begin{array}{ll}
u&=\sin\frac{\pi x}{(1+kt)^{\alpha}}(1+kt)^{-\frac{\alpha}{2}-
\frac{\pi^{2}}{k}}e^{\frac{k^{2}\alpha^{2}x^{2}t}{4(1+kt)}}e^{-\frac{k\alpha x^{2}}{4}}.
\end{array}
\end{equation}

Now let us estimate the norm
\begin{equation}
\label{b10}
\begin{array}{ll}
\|u\|^{2}_{L^{2}(0,l(t))}
&= \int^{(1+kt)^\alpha}_{0}(\sin\frac{\pi x}{(1+kt)^{\alpha}})^{2}{(1+kt)}^{-\alpha}(1+kt)^{-\frac{2\pi^{2}}{k}}
e^{\frac{k^{2}\alpha^{2}x^{2}t}{2(1+kt)}}e^{-\frac{k\alpha x^{2}}{2}}dx\\[3mm]
&\geq(1+kt)^{-\frac{2\pi^{2}}{k}}(1+kt)^{-\alpha}\int^{(1+kt)^\alpha}_{0}(\sin\frac{\pi x}{(1+kt)^{\alpha}})^{2}e^{-\frac{k\alpha x^{2}}{2}}dx.
\end{array}
\end{equation}
Thanks to the following inequality
$$ |\sin\theta|\geq\frac{2}{\pi}|\theta|, \theta\in [-\frac{\pi}{2},\frac{\pi}{2}],$$
we get
\begin{equation}
\label{b11}
\begin{array}{ll}
\|u\|^{2}_{L^{2}(0,l(t))}&\geq\|u\|^{2}_{L^{2}(0,\frac{l(t)}{2})}\\[3mm]
&\geq (1+kt)^{-\frac{2\pi^{2}}{k}}(1+kt)^{-\alpha}\int^{\frac{(1+kt)^\alpha}{2}}_{0}\frac{4}{\pi^{2}}({\frac{\pi x}{(1+kt)^{\alpha}}})^{2}e^{-\frac{k\alpha x^{2}}{2}}dx \\[3mm]
&\geq 4 (1+kt)^{-\frac{2\pi^{2}}{k}}(1+kt)^{-3\alpha}\int^{\frac{(1+kt)^\alpha}{2}}_{0} x^{2}e^{-\frac{k\alpha x^{2}}{2}}dx.
\end{array}
\end{equation}
Set $x=\sqrt{\frac{2}{k\alpha}}y$, the integral in (2.7) turns to be
\begin{equation}
\label{b12}
\begin{array}{ll}
\int^{\frac{(1+kt)^\alpha}{2}}_{0} x^{2}e^{-\frac{k\alpha x^{2}}{2}}dx
&= ({\frac{2}{k\alpha}})^{\frac{3}{2}}\int_{0}^{\sqrt{\frac{k\alpha}{2}}
\frac{(1+kt)^\alpha}{2}}y^{2}e^{-y^{2}}dy\\[3mm]
&=({\frac{2}{k\alpha}})^{\frac{3}{2}}(-\frac{y}{2}e^{-y^{2}}|_{0}^
{\sqrt{\frac{k\alpha}{2}}\frac{(1+kt)^\alpha}{2}}+
\int^{\sqrt{\frac{k\alpha}{2}}\frac{(1+kt)^{\alpha}}{2}}_{0}\frac{e^{-y^{2}}}{2}dy)\\[3mm]
&\geq ({\frac{2}{k\alpha}})^{\frac{3}{2}}(-\frac{1}{8}+\frac{\sqrt{\pi}}{8}) \\[3mm]
&\geq {C_1 },
\end{array}
\end{equation}
where $C_{1}=({\frac{2}{k\alpha}})^{\frac{3}{2}}(-\frac{1}{8}+\frac{\sqrt{\pi}}{8}).$
The inequalities
$$-\frac{y}{2}e^{-y^{2}}|_{y=\sqrt{\frac{k\alpha}{2}}\frac{(1+kt)^\alpha}{2}}\geq -\frac{1}{8}$$
$$\int^{{\sqrt{\frac{k\alpha}{2}}\frac{(1+kt)^{\alpha}}{2}}}_{0}\frac{e^{-y^{2}}}{2}dy\geq \frac{\sqrt{\pi}}{8}$$
hold for $t\geq t_{0}$ due to the facts
$$\lim_{y\rightarrow\infty} ye^{-y^{2}}=0, \int^{\infty}_{0}e^{-y^{2}}dy=\frac{\sqrt{\pi}}{2}.$$
Combining (2.7) and (2.8), we derive that
$$\|u\|_{L^2(0,l(t))}\geq C_{2} (1+kt)^{-\frac{\pi^{2}}{k}}(1+kt)^{-\frac{3\alpha}{2}},$$
where $C_{2}=2\sqrt{C_{1}}.$ 

In the case of $\alpha>\frac{1}{2}$, we can get the similar estimate. Integrating (2.3), we have
\begin{equation}
\label{b7}
\left\{\begin{array}{ll}
  c(x,t)- c(0,t)=-\frac{k\alpha x^{2}}{4(1+kt)},\\[3mm]
 c(x,t)- c(x,0)=-\frac{\pi^{2}}{k(1-2\alpha)}((1+kt)^{(1-2\alpha)}-1)
 -\frac{k\alpha^{2}x^{2}}{4(1+kt)}+\frac{k\alpha^{2}x^{2}}{4}-\frac{\alpha}{2}\ln(1+kt).
\end{array}
\right.
\end{equation}
Taking $c(0,0)=0$, we obtain that
$$c(x,t)=-\frac{\pi^{2}}{(1-2\alpha)k}((1+kt)^{(1-2\alpha)}-1)-\frac{\alpha}{2}\ln(1+kt)-
\frac{k\alpha^{2}x^{2}}{4(1+kt)}+\frac{k\alpha^{2}x^{2}}{4}-\frac{k\alpha x^{2}}{4}.$$
Hence, the solution to (2.1) is
\begin{equation}
\label{b8}
\begin{array}{ll}
u&=\sin\frac{\pi x}{(1+kt)^{\alpha}}e^{-\frac{\pi^{2}}{k(1-2\alpha)}((1+kt)^{(1-2\alpha)}-1)}(1+kt)^
{-\frac{\alpha}{2}}e^{\frac{k^{2}\alpha^{2}x^{2}t}{4(1+kt)}}e^{-\frac{k\alpha x^{2}}{4}}.
\end{array}
\end{equation}
By similar estimate, the norm of the solution turns out to be
\begin{equation}
\label{b10}
\begin{array}{ll}
\|u\|^{2}_{L^{2}(0,l(t))}
&\geq \int^{\frac{(1+kt)^\alpha}{2}}_{0}(\sin\frac{\pi x}{(1+kt)^{\alpha}})^{2}{(1+kt)}^{-\alpha}e^{-\frac{2 \pi^{2}}{k(1-2\alpha)}((1+kt)^{(1-2\alpha)}-1)}e^
{\frac{k^{2}\alpha^{2}x^{2}t}{2(1+kt)}}e^{-\frac{k\alpha x^{2}}{2}}dx\\[3mm]
&\geq e^{-\frac{2 \pi^{2}}{k(1-2\alpha)}((1+kt)^{(1-2\alpha)}-1)}(1+kt)^{-\alpha}\int^{\frac{(1+kt)^\alpha}{2}}_{0}(\sin\frac{\pi x}{(1+kt)^{\alpha}})^{2}e^{-\frac{k\alpha x^{2}}{2}}dx\\[3mm]
&\geq C_{3} e^{-\frac{2 \pi^{2}}{k(1-2\alpha)}((1+kt)^{(1-2\alpha)}-1)}(1+kt)^{-\alpha}.
\end{array}
\end{equation}
We derive that
\begin{equation}
\label{b11}
\begin{array}{ll}
\|u\|_{L^2(0,l(t))}&\geq C_{4} e^{-\frac{ \pi^{2}}{k(1-2\alpha)}((1+kt)^{(1-2\alpha)}-1)}(1+kt)^{-\frac{\alpha}{2}}\\[3mm]
&\geq C_{4} e^{\frac{ \pi^{2}}{k(1-2\alpha)}} (1+kt)^{-\frac{\alpha}{2}}.
\end{array}
\end{equation}

(2)Assume $\alpha=\frac{1}{2}$. Let ${v}_{0}=A\sin(\pi x)e^{-\frac{\alpha x^{2}}{4}},A>0$, then the corresponding solution ${v}$ is only polynomially stable. For any ${u} _{0}\in C(\overline{\Omega}),{u} _{0}>0$, there exists a ${v} _{0} $ such that ${u}_{0}\geq {v}_{0}$. Let $w={u}-{v}$, we have that $w$ is the solution of the following system
\begin{equation}
\label{b13}
\left\{\begin{array}{ll}
w_{t}-w_{xx}=0,
&\mbox{in}\ {Q}_t,\\[3mm]
w(0,t)=w(l(t),t)=0, &\mbox{in}\ (0,+\infty), \\[3mm]
w(x,0)=w_{0}(x)={u}_{0}(x)-{v}_{0}(x)\geq0, &\mbox{in}\ \Omega=(0,1),
\end{array}
\right.
\end{equation}
where ${Q}_t=\{(x,t)| x\in (0,l(t)),t\in(0,+\infty), l(t)=(1+kt)^{\alpha}\}$.
If we reduce the system (2.13) to a variable coefficient parabolic equation in the cylindrical domain, it is clear that the comparison principle holds. Therefore $w\geq0 $ in ${Q}_t $.

Since ${v}\geq0$ in ${{Q}}_t$ , we get ${u}\geq {v}\geq0$ in ${{Q}}_t$. Then we arrive at
$$\|{u}\|_{L^2(0,l(t))}\geq\|{v}\|_{L^2(0,l(t))}\geq C_{2}(1+kt)^{-\frac{\pi^{2}}{k}}(1+kt)^{-\frac{3\alpha}{2}} ,t\geq t_{0}.$$

If $\alpha>\frac{1}{2}$, ${u} _{0}\in C(\overline{\Omega}),{u} _{0}>0$, we have similar estimate.
 
 Hence, we can take $\beta=\frac{\pi^{2}}{k}+\frac{3\alpha}{2}$ for  $\alpha=\frac{1}{2}$ and $\beta=\frac{\alpha}{2}$ for $\alpha>\frac{1}{2}$.
\hfill$\Box$
\begin{corollary}\label{02}
 Let us assume that $u_{0}\in C(\overline{\Omega})$, $u_{0}>0$,$0<\alpha<\frac{1}{2}$. Then
there is some $t_{0}$ such that the solution of system (1.1) satisfies
$$\|u\|_{L^2(0,l(t))}\geq C_{1} (1+kt)^{-\frac{\alpha}{2}}e^{-C_{2}t^{1-2\alpha}} ,t\geq t_{0},$$
where $ C_{1},C_{2}>0$ depending on $k,\alpha, t_{0}$.
\end{corollary}

\noindent{\bf  Proof of Corollary 2.1:}
The proof can be obtained by a simple modification of the first step in the proof of Lemma 2.1. Indeed, integrating (2.3) for $0<\alpha<\frac{1}{2}$, we also have (2.9),(2.10) and (2.11). Namely,
\begin{equation}
\|u\|^{2}_{L^{2}(0,l(t))}\geq  C_{3} e^{-\frac{2 \pi^{2}}{k(1-2\alpha)}((1+kt)^{(1-2\alpha)}-1)}(1+kt)^{-\alpha}.
\end{equation}
Due to $0<1-2\alpha<1$, we get that

\begin{equation}
\label{b11}
\begin{array}{ll}
\|u\|_{L^2(0,l(t))}&\geq C_{4} (1+kt)^{-\frac{\alpha}{2}}e^{-C_{2}(1+kt)^{1-2\alpha}}\\[3mm]
&\geq C'_{4} (1+kt)^{-\frac{\alpha}{2}}e^{-Ct^{1-2\alpha}},t\geq t_{0}.
\end{array}
\end{equation}
Thanks to the inequality
$$ (a+b)^\alpha\leq a^\alpha+b^{\alpha},\quad\mbox{if}\quad a>0,b>0,0<\alpha<1,$$
the last inequality in (2.15) holds.
\hfill$\Box$

\noindent{\bf Proof of Theorem 1.1:}

 Set $x=(1+kt)^{\alpha}y, y\in(0,1), w(y,t)=u((1+kt)^{\alpha}y,t)$. System (1.1) can be rewritten as
\begin{equation}
\label{b3}
\left\{\begin{array}{ll}
 w_{t}-\frac{k\alpha yw_{y}}{1+kt}-\frac{w_{yy}}{(1+kt)^{2\alpha}}=0,
&\mbox{in}\ {Q},\\[3mm]
w(0,t)=w(1,t)=0, &\mbox{in}\ (0,+\infty), \\[3mm]
w(y,0)=u_{0}(y), &\mbox{in}\ \Omega=(0,1),
\end{array}
\right.
\end{equation}
where ${Q}=\{(y,t)| y\in (0,1),t\in(0,+\infty)\}$. System (2.16) is a variable coefficient parabolic equation in the cylindrical domain ${Q}$.

The stability will be shown by the classical energy estimate. Multiplying (2.16) by $w$ and integrating with respect to $y$, we get
\begin{equation}
\label{b4}
 \frac{1}{2}\frac{d}{dt}\int ^{1}_{0}w^{2}dy+ \frac{k\alpha}{2(1+kt)}\int ^{1}_{0}w^{2}dy+\frac{1}{(1+kt)^{2\alpha}}\int ^{1}_{0}w_{y}^{2}dy=0.
\end{equation}
Let $E(t)=\int ^{1}_{0}w^{2}dy$. By Poincar\'{e}'s inequality, we see that
\begin{equation}
\label{b4}
\frac{1}{2}\frac{d}{dt}E(t)+ \frac{k\alpha}{2(1+kt)}E(t)+\frac{C}{(1+kt)^{2\alpha}}E(t)\leq 0
\end{equation}
for a suitable $C>0$.

Using Gronwall inequality, we deduce that for $\alpha=\frac{1}{2}$
\begin{equation}
\label{b1}
 E(t)\leq E(0)(1+kt)^{-(\alpha+\frac{2C}{k})}.
\end{equation}
Consequently,
$$ \|u\|_{L^2(0,l(t))}= E(t)^{ \frac{1}{2}}\leq\|u_{0}\|_{L^2(0,1)}(1+kt)^{-(\frac{\alpha}{2}+\frac{C}{k})}
\leq\|u_{0}\|_{L^2(0,1)}(1+kt)^{-\frac{\alpha}{2}}.$$
Namely, the solution is polynomially stable in the sense of $L^{2}$ norm.

For $\alpha>\frac{1}{2} $, the energy estimate is
\begin{equation}
\label{b4}
 E(t)\leq E(0)(1+kt)^{-\alpha} e^{\frac{2C(1-(1+kt)^{1-2\alpha})}{k(1-2\alpha)}}\leq E(0) (1+kt)^{-\alpha}.
\end{equation}
Thus,$$ \|u\|_{L^2(0,l(t))}\leq\|u_{0}\|_{L^2(0,1)}(1+kt)^{-\frac{\alpha}{2}}.$$

By the estimates above together with Lemma 2.1, we finish the proof of Theorem 1.1.
\hfill$\Box$

\noindent{\bf Proof of Theorem 1.2:}

From (2.18), we have for $0<\alpha<\frac{1}{2}$,
\begin{equation}
\label{b4}
 E(t)\leq E(0)(1+kt)^{-\alpha} e^{\frac{2C(1-(1+kt)^{1-2\alpha})}{k(1-2\alpha)}}\leq E(0) e^{\frac{2C(1-(1+kt)^{1-2\alpha})}{k(1-2\alpha)}}.
\end{equation}
Thanks to $(1+kt)^{1-2\alpha}\geq(kt)^{1-2\alpha}$, we get
$$ \|u\|_{L^2(0,l(t))}\leq\|u_{0}\|_{L^2(0,1)} e^{\frac{C(1-(1+kt)^{1-2\alpha})}{k(1-2\alpha)}}\leq C_{1} e^{-C_{2}t^{1-2\alpha}}.$$
Hence, in the case of $0<\alpha<\frac{1}{2}$, $(1.1)$ is analogously exponentially stable in the sense of $L^{2}$ norm.

Then, in view of Corollary 2.1, we complete the proof of Theorem 1.2.

\hfill$\Box$
\begin{corollary}\label{02}
 Let
 $$l_{1}(t)=(1+kt)^{\alpha_1}, l_{2}(t)=(1+kt)^{\alpha_2} (0<\alpha_{1}<\alpha_{2}),$$
 $$\Omega^{1}_{t}=\{x| 0<x< l_{1}(t)\}, \Omega^{2}_{t}=\{x| 0< x <  l_{2}(t)\},$$
$$u_{0}\in H^{1}_{0}(0,1), u_{0}\geq 0( u_{0}\leq 0).$$
Suppose $L(t)$ is a smooth curve between line $ l_{1}(t)$ and $ l_{2}(t)$, $L(0)=1$. Let $\Omega_{t}=\{x| 0<x< L(t)\}$. If the solutions  corresponding to $ l_{1}(t)$ and $ l_{2}(t)$ are
only polynomially (or analogously exponentially) stable, then the solution of $(1.1)$ with $u_{0}$ and  boundary curve $L(t)$  is only polynomially (or analogously exponentially) stable.
\end{corollary}
\noindent{\bf Proof of Corollary 2.2:}

 Denote by $u_{1}$ the solution corresponding to $u_{0}$ and  $ l_{1}(t)$, and $u_{2}$ the solution corresponding to $u_{0}$ and  $ l_{2}(t)$. Assume $u_{1}$ and $u_{2}$ are
only polynomially stable. Denote the solution corresponding to $u_{0}$ and  $ L(t)$ by $v$. We have that $v|_{\Omega^{1}_{t}\times(0,\infty)}$ is the solution of the following system
\begin{equation}
\label{b1}
\left\{\begin{array}{ll}
v_{t}-\Delta v=0,
&\mbox{in}\ \Omega^{1}_{t}\times(0,\infty),\\[3mm]
v(0,t)=0,v(l_{1}(t),t)=\alpha(t), &\mbox{in}\ (0,+\infty), \\[3mm]
v(x,0)=u_{0}(x)\geq0, &\mbox{in}\ \Omega.
\end{array}
\right.
\end{equation}
 Therefore $\alpha(t)\geq0 $ by the comparison principle for parabolic equation in cylindrical domain and the regularity of the solution in Lemma 1.1. Using the comparison principle again, we get $ v\geq u_{1}\geq0$ in $\Omega^{1}_{t}\times(0,\infty)$. Then we arrive
$$\|v\|_{L^2(\Omega_{t})}\geq\|v\|_{L^2(\Omega^{1}_{t})}\geq\|u_{1}\|_{L^2(\Omega^{1}_{t})}\geq C_{1} (1+kt)^{-\beta} ,t\geq t_{0}.$$
 Similarly,
 $$C_{2}(1+kt)^{-\frac{\alpha_{2}}{2}}\geq\|u_{2}\|_{L^2(\Omega^{2}_{t})}\geq\|v\|_{L^2(\Omega_{t})} ,t\geq t_{0}.$$
The solution $v$ is only polynomially stable.
\hfill$\Box$
\section{Exponential stabilization}
\medskip
In this part, we follow the standard procedure of the backstepping method.

\noindent{\bf  Step1  The stability of the target system}

We choose the stable target system
\begin{equation}
\label{b15}
\left\{\begin{array}{ll}
 w_{t}-w_{xx}+\lambda w=0,
&\mbox{in}\ {Q}_t,\\[3mm]
w(0,t)=0, w(l(t),t)=0, &\mbox{in}\ (0,+\infty), \\[3mm]
w(x,0)=w_{0}(x), &\mbox{in}\ \Omega=(0,1),
\end{array}
\right.
\end{equation}
where ${Q}_t=\{(x,t)| x\in (0,l(t)),t\in(0,+\infty), l(t)=(1+kt)^{\alpha}, k>0, \alpha>0\}, \lambda> 0$ will be determined later.

\begin{lemma}\label{01}
System (3.1) is exponentially stable.
\end{lemma}

We present the proof of Lemma 3.1 which also can be found in [14] as a matter of convenience.

Proof of Lemma 3.1: Multiplying the equation in system $(3.1)$ by $w$, integrating from 0 to $l(t)$ with respect to ${x}$, we get
$$ \frac{1}{2}\frac{d}{dt}\int_{0}^{l(t)}w^{2}dx-ww_{x}|_{0}^{l(t)}+\int_{0}^{l(t)}w^{2}_{x}dx+\lambda\int_{0}^{l(t)}w^{2}dx=0.$$
Taking account into the boundary condition, we have
\begin{equation}
\label{b16}
\begin{array}{ll}
\frac{1}{2}\frac{d}{dt}\int_{0}^{l(t)}w^{2}dx&=-\int_{0}^{l(t)}w^{2}_{x}dx-\lambda\int_{0}^{l(t)}w^{2}dx\\[3mm]
&\leq -\lambda\int_{0}^{l(t)}w^{2}dx.
\end{array}
\end{equation}
We will derive the stability result by Gronwall inequality
$$ \|w\|_{L^{2}(0,l(t))}\leq\|w_{0}\|_{L^{2}(0,1)}e^{-\lambda t}.$$
\medskip
\noindent{\bf  Step2  The equation of the kernel}

We introduce the Volterra transformation as follows,
$$ w(x,t)=u(x,t)+\int_{0}^{x}p(x,y,t)u(y,t)dy. $$
Computing the partial derivatives directly, one gets
\begin{equation}
\label{b17}
w_{x}=u_{x}(x,t)+p(x,x,t)u(x,t)+\int_{0}^{x}p_{x}(x,y,t)u(y,t)dy,
\end{equation}

\begin{equation}
\label{b18}
w_{xx}=u_{xx}+\frac{dp(x,x,t)}{dx}u+p(x,x,t)u_{x}+p_{x}(x,x,t)u+\int_{0}^{x}p_{xx}udy,
\end{equation}

\begin{equation}
\label{b19}
w_{t}=u_{t}+\int_{0}^{x}p(x,y,t)u_{t}(y,t)dy+\int_{0}^{x}p_{t}udy.
\end{equation}
Taking account into the equation in (1.2) and integrating by parts, one has
\begin{equation}\label{b20}
\begin{array}{ll}
w_{t}&=u_{t}(x,t)+\int_{0}^{x}p_{t}udy+\int_{0}^{x}pu_{yy}(y,t)dy\\[3mm]
&=u_{t}(x,t)+\int_{0}^{x}p_{t}udy+\int_{0}^{x}p_{yy}(x,y,t)u(y,t)dy\\[3mm]
&+p(x,x,t)u_{x}(x,t)-p(x,0,t)u_{x}(0,t)
-p_{y}(x,x,t)u(x,t)+p_{y}(x,0,t)u(0,t).
\end{array}
\end{equation}
According to (3.1),(3.4), (3.6), and taking $p(x,0,t)=0$, we derive that
\begin{equation}\label{b21}
\begin{array}{ll}
w_{t}-w_{xx}+\lambda w=-(2\frac{dp(x,x,t)}{dx}-\lambda )u+\int_{0}^{x}(p_{t}-p_{xx}+p_{yy}+\lambda p)udy
\end{array}.
\end{equation}
Therefore, we choose kernel $p(x,y,t)$ defined on $ \mathbb{S}(t)=\{(x,y)|0\leq y\leq x\leq l(t)\}$ satisfying the following system,
\begin{equation}
\label{b19}
\left\{\begin{array}{ll}
p_{t}-p_{xx}+p_{yy}+\lambda p=0,\\[3mm]
p(x,0,t)=0, \\[3mm]
\frac{dp(x,x,t)}{dx}=\frac{\lambda}{2}.
\end{array}
\right.
\end{equation}
The stabilization problem is changed to the existence of the kernel. Meanwhile, the control is the following  feedback of the state by the boundary condition for $w$,
$$ U(t)=-\int_{0}^{l(t)}p(l(t),y,t)u(y,t)dy. $$
\medskip
\noindent{\bf  Step3  The well-posedness and the estimate of the kernel}

In [4], the backstepping method is extended to plants with time-varying coefficients and the explicit expression of the kernel is given. The kernel $p(x,y,t)$ defined on $ \mathbb{D}(t)=\{(x,y)|0\leq y\leq x\leq 1\}$ satisfying(3.8),which is of the following form
\begin{equation}
\label{b18}
p(x,y,t)=-\frac{y}{2}e^{-\lambda t}f(z,t),z=\sqrt{x^{2}-y^{2}},
\end{equation}
where
\begin{equation}
\label{b18}
\left\{\begin{array}{ll}
 f_{t}=f_{zz}+\frac{3f_{z}}{z},\\[3mm]
f_{z}(0,t)=0,f(0,t)=-\lambda e^{\lambda t}:=F(t) .
\end{array}
\right.
\end{equation}
The $C_{z,t}^{2,1}$ solution to this problem is
\begin{equation}
\label{b18}
f(z,t)=\Sigma_{n=0}^{\infty}\frac{1}{n!(n+1)!}(\frac{z}{2})^{2n}F^{(n)}(t).
\end{equation}
Taking account into the form of the kernel in [4], we can obtain the growth order of $p(x,y,t)$. The solution of (3.8) is
\begin{equation}
\label{b22}
\begin{array}{ll}
p(x,y,t)&=-\frac{y}{2}e^{-\lambda t}f(z,t)\\[3mm]
&=-\frac{y}{2}e^{-\lambda t}\Sigma_{n=0}^{\infty}\frac{1}{n!(n+1)!}(\frac{z}{2})^{2n}F^{(n)}(t)\\[3mm]
&=\frac{y}{2}e^{-\lambda t}\Sigma_{n=0}^{\infty}\frac{1}{n!(n+1)!}(\frac{z}{2})^{2n}{\lambda}^{n+1}e^{\lambda t}\\[3mm]
&=\frac{y}{2}\Sigma_{n=0}^{\infty}\frac{1}{n!(n+1)!}(\frac{x^{2}-y^{2}}{4})^{n}{\lambda}^{n+1}\\[3mm]
&=\frac{y}{2l(t)}\Sigma_{n=0}^{\infty}\frac{1}{n!(n+1)!}(\frac{\frac{x^{2}}{l^{2}(t)}-\frac{y^{2}}{l^{2}(t)}}{4})^{n}l^{2n+1}(t){\lambda}^{n+1}.
\end{array}
\end{equation}
The absolute value of the kernel is
\begin{equation}
\label{b22}
\begin{array}{ll}
|p(x,y,t)|
&\leq\frac{1}{2}\Sigma_{n=0}^{\infty}\frac{{\lambda}^{n+1} l^{2n+1}(t)}{4^{n}n!(n+1)!}\\[3mm]
&\leq{\lambda}^{\frac{1}{2}}\Sigma_{n=0}^{\infty}\frac{({\frac{{\lambda}^{\frac{1}{2}} l(t)}{2}})^{n}}{n!}\frac{(\frac{{\lambda}^{\frac{1}{2}} l(t)}{2})^{n+1}}{(n+1)!}\\[3mm]
&\leq{\lambda}^{\frac{1}{2}}\Sigma_{n=0}^{\infty}\frac{({\frac{{\lambda}^{\frac{1}{2}} l(t)}{2}})^{n}}{n!}\Sigma_{n=0}^{\infty}\frac{({\frac{{\lambda}^{\frac{1}{2}} l(t)}{2}})^{n+1}}{(n+1)!}\\[3mm]
&\leq{\lambda}^{\frac{1}{2}}e^{\frac{{\lambda}^{\frac{1}{2}} l(t)}{2}}e^{\frac{{\lambda}^{\frac{1}{2}} l(t)}{2}}\\[3mm]
&\leq{\lambda}^{\frac{1}{2}}e^{{\lambda}^{\frac{1}{2}} l(t)}.
\end{array}
\end{equation}
\medskip
\noindent{\bf  Step4  The invertibility of the transformation}

Let $$u(x,t)=w(x,t)-\int_{0}^{x}q(x,y,t)w(y,t)dy. $$
By similar arguments, the equation for the kernel $q(x,y,t)$ becomes
\begin{equation}
\label{b29}
\left\{\begin{array}{ll}
 q_{t}-q_{xx}+q_{yy}-\lambda q=0,\\[3mm]
q(x,0,t)=0, \\[3mm]
\frac{dq(x,x,t)}{dx}=\frac{\lambda }{2},
\end{array}
\right.
\end{equation}
which is defined on $ \mathbb{S}(t)=\{(x,y)|0\leq y\leq x\leq l(t)\}$.

Proceeding as the analysis of ${k}(x,y,t)$, one can find that the property of the kernel ${q}(x,y,t)$ is similar to the kernel ${k}(x,y,t)$, such as existence, uniqueness and the estimate. Thus, one has
\begin{equation}
\label{b32}
\begin{array}{ll}
|{q}|\leq {\lambda}^{\frac{1}{2}}e^{{\lambda}^{\frac{1}{2}} l(t)}
\end{array}.
\end{equation}

\medskip
\noindent{\bf  Step5  Rapid exponential stabilization of (1.2)}

Now we will show the rapid exponential stabilization of (1.2). It is easy to see by H$\ddot{o}$lder inequality that
\begin{equation}
\label{b34}
\begin{array}{ll}
|\int_{0}^{x}q(x,y,t)w(y,t)dy|
&\leq  {\lambda}^{\frac{1}{2}}e^{{\lambda}^{\frac{1}{2}} l(t)}\int_{0}^{x}|w|dy\\[3mm]
&\leq  {\lambda}^{\frac{1}{2}}e^{{\lambda}^{\frac{1}{2}} l(t)}\|w\|_{L^{2}(0,l(t))}x^{\frac{1}{2}}.
\end{array}
\end{equation}
The estimate for solution of system (1.2) turns out to be
\begin{equation}
\label{b34}
\begin{array}{ll}
\|u\|_{L^{2}(0,l(t))}&\leq\|w\|_{L^{2}(0,l(t))}+\frac{1}{\sqrt{2}}\|w\|_{L^{2}(0,l(t))}{l(t)} {\lambda}^{\frac{1}{2}}e^{{\lambda}^{\frac{1}{2}} l(t)}\\[3mm]
&\leq \|w_{0}\|_{L^{2}(0,1)}e^{-\lambda t}+\frac{1}{\sqrt{2}}{\lambda}^{\frac{1}{2}}l(t)\|w_{0}\|_{L^{2}(0,1)}e^{-\lambda t+{\lambda}^{\frac{1}{2}} (1+kt)^{\alpha}}.
\end{array}
\end{equation}
Since $l(t)$ is unbounded as time $t$ tends to infinity, we need to restrict the growth of $ l(t)$ in order to guarantee the exponential stability of solution. It is readily to verified that the solution is exponentially stable for $0<\alpha\leq1$ because the following inequality
 \begin{equation}
\label{b34}
\begin{array}{ll}
\|u\|_{L^{2}(0,l(t))}&\leq \|w_{0}\|_{L^{2}(0,1)}e^{-\lambda t}+\frac{1}{\sqrt{2}}{\lambda}^{\frac{1}{2}}l(t)\|w_{0}\|_{L^{2}(0,1)}e^{-\lambda t+{\lambda}^{\frac{1}{2}} (1+kt)}\\[3mm]
&\leq e^{-\frac{(\lambda-{\lambda}^{\frac{1}{2}}k) t}{2}}
\end{array}
\end{equation}
 holds for $t\geq t_{0}$ if $t_{0}$ is large enough and $\lambda>k^{2}$.
 Hence, we build a feedback control law to force the solution of the closed-loop system to decay exponentially to zero with arbitrarily prescribed decay rates.

It can be also checked that the feedback control belongs to $H^{1}(0,+\infty)$ for $0<\alpha\leq1.$ From (3.13) and (3.18), we get the following estimate
\begin{equation}
\label{b35}
\begin{array}{ll}
\|U(t)\|^{2}_{L^{2}(0,+\infty)}&=\int^{+\infty}_{0}U^{2}(t)dt\\[3mm]
&=(\int^{+\infty}_{t_0}+\int^{t_{0}}_{0})U^{2}(t)dt\\[3mm]
&\leq(\int^{+\infty}_{t_0}+\int^{t_{0}}_{0})\|u\|^{2}_{L^{2}(0,l(t))}
{({\lambda}^{\frac{1}{2}}e^{{\lambda}^{\frac{1}{2}} l(t)})}^{2}l(t)dt\\[3mm]
&\leq\int^{+\infty}_{t_0} e^{-(\lambda-{\lambda}^{\frac{1}{2}}k) t}{\lambda}e^{2{\lambda}^{\frac{1}{2}} (1+kt)}(1+kt)dt+\int^{t_{0}}_{0} \|u\|^{2}_{L^{2}(0,l(t))}{\lambda}e^{2{\lambda}^{\frac{1}{2}} l(t)}l(t)dt\\[3mm]
&\leq\int^{+\infty}_{t_0} e^{-(\lambda-3{\lambda}^{\frac{1}{2}}k) t}{\lambda}e^{2{\lambda}^{\frac{1}{2}}}(1+kt)dt+\int^{t_{0}}_{0} \|u\|^{2}_{L^{2}(0,l(t))}{\lambda}e^{2{\lambda}^{\frac{1}{2}} l(t)}l(t)dt\\[3mm]
&< \infty,
\end{array}
\end{equation}
provided $\lambda>9k^{2}$, which means that $U(t)\in L^{2}(0,+\infty). $

On the other hand,
\begin{equation}
\label{b36}
\begin{array}{ll}
U'(t)&=-l'(t)p(l(t),l(t),t)u(l(t),t)-\int_{0}^{l(t)}p_{x}(l(t),y,t)l'(t)u(y,t)dy\\[3mm]
&-\int_{0}^{l(t)}p_{t}(l(t),y,t)u(y,t)dy-\int_{0}^{l(t)}p(l(t),y,t)u_{t}(y,t)dy.
\end{array}
\end{equation}

The first term belongs to $L^{2}(0,+\infty)$ provided $\lambda>25k^{2}$due to
\begin{equation}
\label{b37}
\begin{array}{ll}
\int_{t_0}^{+\infty}(l'(t)p(l(t),l(t),t)u(l(t),t))^{2}dt
&=\int_{t_0}^{+\infty}(l'(t)p(l(t),l(t),t)U(t))^{2}dt\\[3mm]
&\leq\int_{t_0}^{+\infty}{k^{2}\alpha}^{2}(1+kt)^{2\alpha-2}({\lambda}^{\frac{1}{2}}e^{{\lambda}^{\frac{1}{2}} l(t)})^{4}l(t) e^{-(\lambda-k{\lambda}^{\frac{1}{2}}) t}dt\\[3mm]
&\leq\int_{t_0}^{+\infty}C(k,\alpha,\lambda)(1+kt)^{3\alpha-2}e^{-(\lambda-5{k\lambda}^{\frac{1}{2}}) t}dt\\[3mm]
&<\infty.
\end{array}
\end{equation}

Using the expression of the kernel, we can calculate the derivatives with respect to $x$ and $t,$
\begin{equation}
\label{b37}
\begin{array}{ll}
p_{x}(x,y,t)&=-\frac{y}{2}e^{-\lambda t}f_{z}(z,t)\frac{x}{\sqrt{x^{2}-y^{2}}}\\[3mm]
&=-\frac{xy}{2z}e^{-\lambda t}\Sigma_{n=1}^{\infty}\frac{n}{n!(n+1)!}(\frac{z}{2})^{2n-1}F^{(n)}(t),
\end{array}
\end{equation}

\begin{equation}
\label{b37}
\begin{array}{ll}
p_{t}(x,y,t)=\frac{y}{2}\lambda e^{-\lambda t}f(z,t)-\frac{y}{2} e^{-\lambda t}f_{t}(z,t).
\end{array}
\end{equation}
The absolute value of $p_{t}(x,y,t)$ is
\begin{equation}
\label{b37}
\begin{array}{ll}
|p_{t}(x,y,t)|&\leq\lambda|p|+|\frac{y}{2}e^{-\lambda t}\Sigma_{n=0}^{\infty}\frac{1}{n!(n+1)!}
(\frac{z}{2})^{2n}\lambda^{n+2}e^{\lambda t}|\\[3mm]
&\leq2\lambda|p|.
\end{array}
\end{equation}
The absolute value of $p_{x}(x,y,t)$ is
\begin{equation}
\label{b37}
\begin{array}{ll}
|p_{x}|&\leq
{\frac{xy}{4}}e^{-\lambda t}\Sigma_{n=1}^{\infty}\frac{1}{n!n!}(\frac{z}{2})^{2n-2}F^{(n)}(t)\\[3mm]
&\leq x^{2}e^{-\lambda t}\Sigma_{n=1}^{\infty}\frac{1}{n!n!}
\frac{(x^{2}-y^{2})^{n-1}}{4^{n}}{\lambda}^{n+1}e^{\lambda t}\\[3mm]
&\leq
(\frac{x}{l(t)})^{2}\Sigma_{n=1}^{\infty}\frac{1}{n!n!4^{n}}
(\frac{x^{2}-y^{2}}{l^{2}(t)})^{n-1}l^{2n}(t){\lambda}^{n+1}\\[3mm]
&\leq\Sigma_{n=1}^{\infty}\frac{1}{n!n!4^{n}}
l^{2n}(t){\lambda}^{n+1}\\[3mm]
&\leq\lambda e^{\lambda^{\frac{1}{2}}l(t)}.
\end{array}
\end{equation}
The absolute value of the second term on the right hand of (3.20) is
\begin{equation}
\label{b37}
\begin{array}{ll}
|\int_{0}^{l(t)}p_{x}(l(t),y,t)l'(t)u(y,t)dy|
&\leq\lambda e^{\lambda^{\frac{1}{2}}l(t)}l'(t)\|u\|_{L^{2}(0,l(t))}\sqrt{l(t)}\\[3mm]
&\leq C(k,\alpha,\lambda)(1+kt)^{\frac{3\alpha}{2}-1} e^{-\frac{(\lambda-3k{\lambda}^{\frac{1}{2}}) t}{2}}.
\end{array}
\end{equation}
When $\lambda>9k^{2}$, the $L^{2}$ norm estimate of the second term becomes
\begin{equation}
\label{b37}
\begin{array}{ll}
\int_{0}^{+\infty}(\int_{0}^{l(t)}p_{x}(l(t),y,t)l'(t)u(y,t)dy)^{2}dt
&\leq\int_{0}^{+\infty}C'(k,\alpha,\lambda)(1+kt)^{3\alpha-2} e^{-(\lambda-3k{\lambda}^{\frac{1}{2}}) t}dt\\[3mm]
&<\infty.
\end{array}
\end{equation}
Similarly, the third term in (3.20) belongs to $L^{2}(0,+\infty)$ by H$\ddot{o}$lder inequality.

Now let us deal with the fourth term in (3.20). According to the inverse transformation $$u(x,t)=w(x,t)-\int_{0}^{x}q(x,y,t)w(y,t)dy,$$
the derivative of $u$ with respect to $t$ is
$$u_{t}(x,t)=w_{t}(x,t)-\int_{0}^{x}q_{t}(x,y,t)w(y,t)dy
-\int_{0}^{x}q(x,y,t)w_{t}(y,t)dy.
$$
The $L^{2}$ norm of $u_{t}$ satisfies
\begin{equation}
\label{b37}
\begin{array}{ll}
\|u_{t}\|_{L^{2}(0,l(t))}
&\leq\|w_{t}\|_{L^{2}(0,l(t))}+\bar{q}_{t}\sqrt{l(t)}\|w\|_{L^{2}(0,l(t))}
+\bar{q}\sqrt{l(t)}\|w_{t}\|_{L^{2}(0,l(t))}\\[3mm]
&\leq\bar{q}_{t}\sqrt{l(t)}\|w\|_{L^{2}(0,l(t))}
+2\bar{q}\sqrt{l(t)}\|w_{t}\|_{L^{2}(0,l(t))}\\[3mm]
&\leq2(\lambda+1)\bar{q}\sqrt{l(t)}(\|w\|_{L^{2}(0,l(t))}
+\|w_{t}\|_{L^{2}(0,l(t))}),
\end{array}
\end{equation}
where $\bar{q}_{t}$ denotes the upper bound of $|{q}_{t}|$, $\bar{q}=\bar{p}={\lambda}^{\frac{1}{2}}e^{{\lambda}^{\frac{1}{2}} l(t)}$ .
The last inequality in (3.28) is the consequence of $\bar{q}_{t}\leq2\lambda\bar{q}$.
In order to estimate (3.28), let us introduce the change of variable in system (3.1)
$$w=\tilde{w}e^{-\lambda t},$$
then $\tilde{w}$ is the solution to
\begin{equation}
\label{b15}
\left\{\begin{array}{ll}
 \tilde{w}_{t}-\tilde{w}_{xx}=0,
&\mbox{in}\ {Q}_t,\\[3mm]
\tilde{w}(0,t)=0, \tilde{w}(l(t),t)=0, &\mbox{in}\ (0,+\infty), \\[3mm]
\tilde{w}(x,0)=w_{0}(x), &\mbox{in}\ \Omega=(0,1).
\end{array}
\right.
\end{equation}
If $w_{0}\in H^{1}_{0}(0,1)$, system (3.29) has a unique strong solution. And we have the following estimates by Lemma 4.2 and the change of variable,
\begin{equation}
\|\tilde{w}\|_{L^{2}(0,l(t))}+\|\tilde{w}_{t}\|_{L^{2}(Q_{t})}\leq C\|{w}_{0}\|_{H^{1}_{0}(0,1)},
\end{equation}
\begin{equation}
\|{w}\|_{L^{2}(0,l(t))}\leq Ce^{-\lambda t}\|{w}_{0}\|_{H^{1}_{0}(0,1)},
\end{equation}
\begin{equation}
\|{w}_{t}\|_{L^{2}(0,l(t))}\leq e^{-\lambda t}\|{\tilde{w}}_{t}\|_{L^{2}(0,l(t))}+\lambda e^{-\lambda t}\|{\tilde{w}}\|_{L^{2}(0,l(t))}.
\end{equation}
Now we will estimate the last term in (3.20). By (3.28) and H$\ddot{o}$lder inequality, we get
\begin{equation}
\label{b37}
\begin{array}{ll}
&\int_{0}^{+\infty}(\int_{0}^{l(t)}p(l(t),y,t)u_{t}(y,t)dy)^{2}dt\\[3mm]
&\leq\int_{0}^{+\infty}\bar{p}l(t)\|u_{t}\|^{2}_{L^{2}(0,l(t))}dt\\[3mm]
&\leq\int_{0}^{+\infty}\bar{p}l(t)8(\lambda+1)^{2}\bar{q}^{2}l(t)(\|w\|^{2}_{L^{2}(0,l(t))}
+\|w_{t}\|^{2}_{L^{2}(0,l(t))})dt
\\[3mm]
&\leq C(\lambda)\int_{0}^{+\infty}e^{3{\lambda}^{\frac{1}{2}} l(t)}l^{2}(t)(\|w\|^{2}_{L^{2}(0,l(t))}
+\|w_{t}\|^{2}_{L^{2}(0,l(t))})dt.
\end{array}
\end{equation}
Combining (3.30), (3.31) with (3.32), we obtain that
\begin{equation}
\label{b37}
\begin{array}{ll}
&\int_{0}^{+\infty}(\int_{0}^{l(t)}p(l(t),y,t)u_{t}(y,t)dy)^{2}dt\\[3mm]
&\leq C(\lambda)\int_{0}^{+\infty}e^{3{\lambda}^{\frac{1}{2}} l(t)}l^{2}(t)(\|w\|^{2}_{L^{2}(0,l(t))}
+2e^{-2\lambda t}\|{\tilde{w}}_{t}\|^{2}_{L^{2}(0,l(t))}+2\lambda ^{2}e^{-2\lambda t}\|{\tilde{w}}\|^{2}_{L^{2}(0,l(t))})dt\\[3mm]
&\leq C'(\lambda)\int_{0}^{+\infty}e^{3{\lambda}^{\frac{1}{2}} l(t)-\lambda t}l^{2}(t)\|{w}_{0}\|^{2}_{H^{1}_{0}(0,1)}dt
+C'(\lambda)\int_{0}^{t_{0}}e^{3{\lambda}^{\frac{1}{2}} l(t)}l^{2}(t)e^{-2\lambda t}\|{\tilde{w}}_{t}\|^{2}_{L^{2}(0,l(t))}dt\\[3mm]
&+C'(\lambda)\int_{t_{0}}^{+\infty}\|{\tilde{w}}_{t}\|^{2}_{L^{2}(0,l(t))}dt\\[3mm]
&<+\infty.
\end{array}
\end{equation}

Thus, we derive that $U(t)\in H^{1}(0,+\infty)$ if $w_{0}\in H^{1}_{0}(0,1)$. If $w_{0}\in L^{2}(0,1)$, the solution of system (3.29) will satisfies $w(\cdot,T)\in H^{1}_{0}(0,l(T))$ at some time $T$ by the regularity of heat equation.

Similarly, one can also prove $\sqrt{l(t)}U'(t)\in L^{2}(0,+\infty)$ for large enough $\lambda$.

In conclusion, we should take $\lambda>25k^{2}$ in order to guarantee the exponential stabilizability, $U(t)\in H^{1}(0,+\infty)$ and $\sqrt{l(t)}U'(t)\in L^{2}(0,+\infty)$.
\begin{remark}\label{03}
The well-posedness of solution to system (1.2) with the feedback control can be see from the well-posedness of solutions to system (3.1) and system (3.8).
\end{remark}
\begin{remark}\label{03}
When the growth order of $l(t)$ is less than 1, such as  $l(t)=1+ln(1+t),l(t)=1+sint$, the linear feedback control does work from (3.17).
\end{remark}
\begin{remark}\label{03}
 Since the PDE on the non-cylindrical domain can be converted to the equation with time-dependent coefficients, the results of this paper may be extended to some parabolic equations with time-dependent coefficients, for which the explicit expression of the kernel can not be obtained. These are problems which we will consider next.
\end{remark}
\hfill$\Box$
\section{Appendix}
Let us start with the following system
\begin{equation}
\label{b1}
\left\{\begin{array}{ll}
 u_{t}-u_{xx}=f(x,t),
&\mbox{in}\ {{Q}}_T,\\[3mm]
u(0,t)=0,u(l(t),t)=0, &\mbox{in}\ (0,T), \\[3mm]
u(x,0)=u_{0}(x), &\mbox{in}\ \Omega=(0,1),
\end{array}
\right.
\end{equation}
where ${{Q}}_T=\{(x,t)| x\in (0,l(t))=\Omega_{t},t\in(0,T), l(t)=(1+kt)^{\alpha},k>0,\alpha>0\}$.
\begin{lemma}([16])
If $u_{0}\in L^{2}(\Omega)$, $f\in L^{2}(0,T;H^{-1}(\Omega_{t})) $, there exists a unique weak solution of (4.1) in the following space
$$u\in C([0,T];L^2(\Omega_{t}))\cap  L^2(0,T;H^1_{0}(\Omega_{t})).$$
Moreover, there exists a positive constant C (independent with ${Q}_{T},u_{0},f$) such that
\begin{equation}\|u\|^{2}_{L^{\infty}(0,T;L^2(\Omega_{t}))}+ \|u\|^{2}_{L^{2}(0,T;H^{1}_{0}(\Omega_t))}\leq C[\|u_{0}\|^{2}_{L^2(\Omega)}+\|f\|^{2}_{L^2(0,T;H^{-1}(\Omega_t))}].
 \end{equation}
\end{lemma}

\begin{lemma}
If $u_{0}\in H^{1}_{0}(\Omega)$, $f\in L^{2}(0,T;L^{2}(\Omega_{t})) $, problem (4.1) has a unique strong solution $u$,
 $$u\in C([0,T];H^1_{0}(\Omega_{t}))\cap  L^2(0,T;H^2\cap H^1_{0}(\Omega_{t}))\cap H^1(0,T;L^2(\Omega_{t})).$$ Moreover, there exists a positive constant C (independent with ${Q}_{T},u_{0},f$) such that
\begin{equation}\|u\|^{2}_{L^{\infty}(0,T;H^{1}_{0}({\Omega}_t))}+ \|u\|^{2}_{L^{2}(0,T;H^{2}({\Omega}_t))}
 \leq C[\|u_{0}\|^{2}_{H^{1}_{0}(\Omega)}+\|f\|^{2}_{L^2(0,T;L^{2}(\Omega_{t}))}].
 \end{equation}
\end{lemma}
\noindent{\bf Proof of Lemma 4.2:}
 The authors ([16]) deduced the energy estimate
 $$ \|u\|^{2}_{L^{\infty}(0,T;H^{1}_{0}({\Omega}_t))}+ \|u\|^{2}_{L^{2}(0,T;H^{2}({\Omega}_t))}
 \leq C[\|u_{0}\|^{2}_{H^{1}_{0}(\Omega)}+\|f\|^{2}_{L^2(Q_{T})}],$$
with $C>0$ dependent of ${Q}_{T}$.

If $f\in L^{2}(0,T;L^{2}(\Omega_{t}))$, a slight variation in the argument allows also to get
\begin{equation}\|u\|^{2}_{L^{\infty}(0,T;H^{1}_{0}({\Omega}_t))}+ \|u\|^{2}_{L^{2}(0,T;H^{2}({\Omega}_t))}
 \leq C[\|u_{0}\|^{2}_{H^{1}_{0}(\Omega)}+\|f\|^{2}_{L^2(0,T;L^{2}(\Omega_{t}))}],
 \end{equation}
with $C>0$ independent of ${Q}_{T}$.
Indeed, instead of using Gronwall inequality to ((14) in [16])
\begin{equation}
\frac{d}{dt}\int_{{\Omega}_t}|\nabla u|^{2}dx\leq-\frac{1}{2}\int_{\Omega_t}|\Delta u|^{2}dx+\| f\|_{L^2(\Omega_{t})}\|\Delta u\|_{L^2(\Omega_{t})}+c\|\nabla u\|^{2}_{L^2(\Omega_{t})},
\end{equation}
we integrate (4.5) with respect to $t$ from 0 to $t$
\begin{equation}
\label{b37}
\begin{array}{ll}
&\int_{{\Omega}_t}|\nabla u|^{2}dx+\frac{1}{2}\int^{t}_{0}\int_{\Omega_t}|\Delta u|^{2}dxdt\\[3mm]
&\leq\int_{{\Omega}}|\nabla u_{0}|^{2}dx+\int^{t}_{0}\| f\|_{L^2(\Omega_{t})}\| \Delta u\|_{L^2(\Omega_{t})}dt+c\int^{t}_{0}\|\nabla u\|^{2}_{L^2(\Omega_{t})}dt.\\[3mm]
\end{array}
\end{equation}
Using the Young inequality, we have
\begin{equation}
\label{b37}
\begin{array}{ll}
&\int_{{\Omega}_t}|\nabla u|^{2}dx+\frac{1}{2}\int^{t}_{0}\int_{\Omega_t}|\Delta u|^{2}dxdt\\[3mm]
&\leq\|\nabla u_{0}\|^{2}_{H^{1}_{0}(\Omega)}+c\|f\|^{2}_{L^2(0,T;L^{2}(\Omega_{t}))}+\frac{1}{4}\int^{t}_{0}\|\Delta u\|^{2}_{L^2(\Omega_{t})}dt+ c\int^{t}_{0}\|\nabla u\|^{2}_{L^2(\Omega_{t})}dt.\\[3mm]
\end{array}
\end{equation}
Then, in view of (4.2) we have
\begin{equation}
\label{b37}
\begin{array}{ll}
\int_{{\Omega}_t}|\nabla u|^{2}dx+\frac{1}{4}\int^{t}_{0}\int_{\Omega_t}|\Delta u|^{2}dxdt\leq C[\|\nabla u_{0}\|^{2}_{H^{1}_{0}(\Omega)}+\|f\|_{L^2(0,T;L^{2}(\Omega_{t}))}].
\end{array}
\end{equation}
for a suitable $C>0$.
\hfill$\Box$

Due to the uniform estimate and non blow-up property of the solution, we get that the existence domain of the solution is $(0,+\infty)$. Lemma 1.1 follows from Lemma 4.1 and Lemma 4.2 directly.
And we also obtain the well-posedness of the following system
\begin{equation}
\label{b1}
\left\{\begin{array}{ll}
 u_{t}-u_{xx}=f(x,t),
&\mbox{in}\ {Q}_t,\\[3mm]
u(0,t)=0,u(l(t),t)=0, &\mbox{in}\ (0,+\infty), \\[3mm]
u(x,0)=u_{0}(x), &\mbox{in}\ \Omega=(0,1),
\end{array}
\right.
\end{equation}
where ${Q}_t=\{(x,t)| x\in (0,l(t))=\Omega_{t},t\in(0,+\infty), l(t)=(1+kt)^{\alpha},k>0,\alpha>0\}$.
 \begin{lemma}
If $u_{0}\in L^{2}(\Omega)$, $f\in L^{2}(0,+\infty;H^{-1}(\Omega_{t})) $, there exists a unique weak solution of (4.9) in the following space
$$u\in C([0,+\infty);L^2(\Omega_{t}))\cap  L^2(0,+\infty;H^1_{0}(\Omega_{t})).$$
Moreover, there exists a positive constant C (independent with ${Q}_{t},u_{0},f$) such that
\begin{equation}\|u\|^{2}_{L^{\infty}(0,+\infty;L^2(\Omega_{t}))}+ \|u\|^{2}_{L^{2}(0,+\infty;H^{1}_{0}(\Omega_t))}\leq C[\|u_{0}\|^{2}_{L^2(\Omega)}+\|f\|^{2}_{L^2(0,+\infty;H^{-1}(\Omega_t))}].
 \end{equation}

If $u_{0}\in H^{1}_{0}(\Omega)$,$f\in L^{2}(0,+\infty;L^{2}(\Omega_{t})) $, problem (4.9) has a unique strong solution
 $$u\in C([0,+\infty);H^1_{0}(\Omega_{t}))\cap  L^2(0,+\infty;H^2\cap H^1_{0}(\Omega_{t}))\cap H^1(0,+\infty;L^2(\Omega_{t})).$$Moreover, there exists a positive constant C (independent with $Q_{t},u_{0},f$) such that
\begin{equation}\|u\|^{2}_{L^{\infty}(0,+\infty;H^{1}_{0}({\Omega}_t))}+ \|u\|^{2}_{L^{2}(0,+\infty;H^{2}({\Omega}_t))}
 \leq C[\|u_{0}\|^{2}_{H^{1}_{0}(\Omega)}+\|f\|^{2}_{L^2(0,+\infty;L^{2}(\Omega_{t}))}].
 \end{equation}
\end{lemma}

\noindent{\bf Proof of Lemma 1.2:} System (1.2) can be converted to system (4.9) by the transformation
$$ \hat{u}=u-\frac{x}{l(t)}U(t).$$ And $\hat{u}$ is the solution of
\begin{equation}
\label{b1}
\left\{\begin{array}{ll}
\hat{ u}_{t}-\hat{u}_{xx}=f(x,t),
&\mbox{in}\ {Q}_t,\\[3mm]
\hat{u}(0,t)=0,\hat{u}(l(t),t)=0, &\mbox{in}\ (0,+\infty), \\[3mm]
\hat{u}(x,0)=\hat{u}_{0}(x), &\mbox{in}\ \Omega=(0,1),
\end{array}
\right.
\end{equation}
where $f(x,t)=-\frac{xU'(t)}{l(t)}+\frac{xl'(t)U(t)}{l^{2}(t)}$. It is easy to check that $f(x,t)\in L^{2}(0,+\infty;L^{2}(\Omega_{t}))$ provided that $U(t)\in H^{1}(0,+\infty)$,$\sqrt{l(t)}U'(t)\in L^{2}(0,+\infty)$ and $\alpha\in (0,2]$. So we prove Lemma 1.2.

\end{document}